\documentclass[11pt]{article}
\usepackage{amssymb}%,psfig}
\usepackage{amsthm}
\usepackage[centertags]{amsmath}
\usepackage{amsfonts}
\usepackage{amssymb}
\usepackage{amsthm}
\usepackage{epsfig}
\usepackage{setspace}
\usepackage{ae} % Umlauts
\usepackage{eucal} % I think standard calligraphics are horrible
\def\citet{\cite}

\newtheorem{theorem}{Theorem}%[section]
\newtheorem{lemma}{Lemma}%[section]
\newtheorem{corollary}{Corollary}%[section]

\newtheorem{remark}{Remark}%[section]

\newcounter{thanksnum}
\def\thanksnumber#1
{\setcounter{thanksnum}{\value{footnote}}\setcounter{footnote}{#1}%
                     \addtocounter{footnote}{-1}\footnotemark
                     \setcounter{footnote}{\value{thanksnum}}}
\def\newtheoremz#1{\@ifnextchar[{\@othmz{#1}}{\@nthmz{#1}}}

\def\@nthmz#1#2{%
\@ifnextchar[{\@xnthmz{#1}{#2}}{\@ynthmz{#1}{#2}}}

\def\@xnthmz#1#2[#3]{\expandafter\@ifdefinable\csname #1\endcsname
{\@definecounter{#1}\@addtoreset{#1}{#3}%
\expandafter\xdef\csname the#1\endcsname{\expandafter\noexpand
  \csname the#3\endcsname \@thmcountersepz \@thmcounterz{#1}}%
\global\@namedef{#1}{\@thmz{#1}{#2}}\global\@namedef{end#1}{\@endtheoremz}}}

\def\@ynthmz#1#2{\expandafter\@ifdefinable\csname #1\endcsname
{\@definecounter{#1}%
\expandafter\xdef\csname the#1\endcsname{\@thmcounterz{#1}}%
\global\@namedef{#1}{\@thm{#1}{#2}}\global\@namedef{end#1}{\@endtheoremz}}}

\def\@othmz#1[#2]#3{\expandafter\@ifdefinable\csname #1\endcsname
  {\global\@namedef{the#1}{\@nameuse{the#2}}%
\global\@namedef{#1}{\@thmz{#2}{#3}}%
\global\@namedef{end#1}{\@endtheoremz}}}

\def\@thmz#1#2{\refstepcounter
    {#1}\@ifnextchar[{\@ythmz{#1}{#2}}{\@xthmz{#1}{#2}}}

\def\@xthmz#1#2{\@begintheoremz{#2}{\csname the#1\endcsname}\ignorespaces}
\def\@ythmz#1#2[#3]{\@opargbegintheoremz{#2}{\csname
       the#1\endcsname}{#3}\ignorespaces}

\def\@thmcounterz#1{\noexpand\arabic{#1}}
\def\@thmcountersepz{.}
\def\@begintheoremz#1#2{ \trivlist \item[\hskip \labelsep{\bf #1\ #2}]}
\def\@opargbegintheoremz#1#2#3{ \trivlist
      \item[\hskip \labelsep{\bf #1\ #2\ (#3)}]}
\def\@endtheoremz{\endtrivlist}

\def\d{\delta}

\def\O{\Omega}

\def\F{{\cal F}}
\def\w{\widehat}

\def\R{{\bf R}}
\def\E{{\bf E}}
\def\P{{\bf P}}

\def\H{{\cal H}}

\def\b{\beta}
\def\s{\delta}
\def\g{\gamma}

\def\t{\theta}
\def\oo{\bar}
\def\s{\sigma}
\def\G{\Gamma}
\def\GG{{\cal G}}

\def\A{{\cal A}}

\newcommand{\be}{\begin{equation}}
\newcommand{\ee}{\end{equation}}
\newcommand{\bd}{\begin{displaymath}}
\newcommand{\ed}{\end{displaymath}}
\newcommand{\ba}{\begin{array}{ll}}
\newcommand{\ea}{\end{array}}
\newcommand{\baa}{\begin{eqnarray}}
\newcommand{\eaa}{\end{eqnarray}}
\newcommand{\baaa}{\begin{eqnarray*}}
\newcommand{\eaaa}{\end{eqnarray*}}   \font\sm=cmr10
\def\H{{\cal H}}

\def\Ts{s}
\def\DD{{\cal D}}
\title{A smooth component of the fractional Brownian motion
and optimal portfolio  selection}
\author{
Nikolai Dokuchaev\\
 {\sm Department of Mathematics \& Statistics, Curtin
University}\\
{\sm  email: N.Dokuchaev@curtin.edu.au } }
\begin{document}
\maketitle
\begin{abstract} We consider fractional Brownian motion with the Hurst parameters from $(1/2,1)$.
We found that the increment of a fractional Brownian motion can be  represented as the sum of a two independent Gaussian
processes  one of which is smooth in the sense that it is differentiable  in mean square.
We  consider fractional Brownian motion and
stochastic integrals generated by the Riemann  sums. As an example of applications,
this results is used to find an optimal pre-programmed strategy in the mean-variance setting for a Bachelier type market model driven by a fractional Brownian
motion.
\end{abstract}
%see(?)http://www.artofproblemsolving.com/LaTeX/AoPS_L_GuidePack.php
 \section{Introduction}
 In this short note, we  consider fractional Brownian motion and
stochastic integrals generated by the Riemann  sums.
We found that the increment of a fractional Brownian motion $B_H(t)$ with the zero mean can be  represented as the sum of a two independent Gaussian
processes  one of  which is smooth in the sense that it is differentiable  in mean square sense,
with the derivative that is square integrable on the finite time intervals. Similarly to the drift part of the diffusion processes, expectations of the
integrals by
this process are non-zero for the processes adapted to it.
 This process can be considered as an analog of the drift. It has to be noted that  the term "drift" is usually applied to
 $\mu$ presented  for the process $\mu t+B_H(t)$; see
\cite{CNS},\cite{Mu},\cite{Es}, where  and estimation of $\mu$ was studied.
In  \cite{Mu}, the term "drift" was also used for a representation for $B_H$ after linear integral transformation and random time change via a standard Brownian motion process with constant in time drift.
 Our representation is for $B_H$ itself, i.e, without an integral transformation.

We also present  an example of applications for portfolio selection.
Statistical properties of fractional Brownian motion are widely used
for financial modelling; see, e.g., the review in \cite{B11} and \cite{H}.
Optimal portfolio selection also have been studied; see, e.g., \cite{Bi}, \cite{J},\cite{H},\cite{H2}, \cite{CS}, as well as more general optimal
stochastic control problems \cite{Zhou}.
It is known that the fractional market with the Hurst parameter $H>1/2$ allows arbitrage (see, e.g., \cite{B11}).
The arbitrage opportunities can be eliminated  by inclusion of  transaction costs in the model \cite{G}.
or by additional restrictions on the strategies such as in \cite{C}.  We consider pre-programmed strategies; this also helps to exclude arbitrage.
We study the mean-variance setting for a Bachelier type market model driven by a fractional Brownian
motion. This linear quadratic optimization setting, or so-called mean-variance portfolio, was
   first introduced for single period models by Markowitz (see \cite{Mar}) and extended later on multi-period and continuous time problems
   (see, e.g., \cite {L2} and \cite{D10} for multi-period problems and \cite{L4} for continuous time).

\section{The  main result}
 We are
given a standard probability space $(\Omega,\F,\P)$, where $\Omega$ is
a set of elementary events, $\F$ is a complete $\s$-algebra of
events, and $\P$ is a probability measure.

 We  assume that $\{B_H(t)\}_{t\in\R}$ is a fractional Brownian motion such that $B_H(0)=0$
 with the Hurst parameter $H\in  (1/2,1)$ defined as described in \cite{Ma,GN}
 such that
 \baa
&& B_H(t)-B_H(s)\nonumber\\&&=c_H\int_s^t(t-q)^{H-1/2}dB(q)+c_H  \int_{-\infty}^s\left[(t-q)^{H-1/2}-(s-q)^{H-1/2}\right]dB(q),
 \label{BBH}
 \eaa
 where $t>s$,
 $c_H=\sqrt{2H\G (3/2-H)/[\G(1/2+H)\G(2-2H)]}$ and  $\G$ is the gamma function. Here
$\{B(t)\}_{t\in\R}$ is  standard Brownian motion such that $B(0)=0$.
Let $\{\GG_t\}$ be the filtration generated by the process $B(t)$.

By (\ref{BBH}),
\baaa
 B_H(t)-B_H(s)=W_H(t)+R_H(t),
 \eaaa
where
 \baaa
 W_H(t)=c_H\int_{s}^t(t-q)^{H-1/2}dB(q),\qquad
R_H(t)=c_H  \int_{-\infty}^{s}f(t,q)dB(q),
\label{BBHa}\eaaa
and where $f(t,q)=(t-q)^{H-1/2}-(s-q)^{H-1/2}$.

Let $s\ge 0$ and $T>s$ be fixed.
For $\tau\in [s,T]$ and  $g\in L_2(s,T)$, set \baaa
G_H(\tau,s,T,g)=c_H(H-1/2)\int_{\tau}^T(t-\tau)^{H-3/2}g(t)dt.
\eaaa
\begin{theorem}\label{ThM}
The processes $W_H(t)$ and $R_H(t)$, where $t>s$, are independent Gaussian $\{\GG_t\}$-adapted processes with zero mean and such that
 the following holds.
 \begin{enumerate}
 \item
 $W_H(t)$ is independent on $\GG_s$ for all $t>s$ and has an It\^o's differential in $t$ in the following sense:  for any $T>s$,
 there exists a function $h(\cdot,s,T)\in L_2(s,T)$ such that 
 \baaa
 \int_{s}^T\g(t)dW_H(t)=\int_{s}^TG_H(\tau,s,T,\g)d B(\tau)
\eaaa
for any $\g\in L_2(\O,\GG_s,\P,L_2(s,T))$.
\item
 $R_H(t)$ is $\GG_s$-measurable for all $t>s$ and differentiable in $t>s$ in mean square sense. More precisely,
 there exits a process   $\DD R_H$ such that
 \subitem(a)  $\DD R_H(t)$ is $\GG_s$-measurable for all $t>s$;
\subitem(b) for any $t>s$, \baaa \E \DD R_H(t)^2=c_H^2\frac{H-1/2}{2}(t-\Ts)^{2H-2},\qquad
\E\int_{s}^t \DD R_H(q)^2dq<+\infty;
\eaaa
\subitem(c)
for any $t>s$,  \baa
\lim_{\d\to 0} \E\left|\frac{R_H(t+\d)-R_H(t)}{\d} -\DD R_H(t)\right|=0.\label{deflim}
\eaa
\end{enumerate}
\end{theorem}
 \begin{corollary}\label{corr} Let $\g(t)$ be a process such that $\g(t)$ is $\GG_s$-measurable for all $t$,
and that $\E\int_s^T\g(t)^2dt<+\infty$ for any $T>s$.
Then the following holds.
\begin{enumerate}
\item[(i)]
\baaa
\int_{s}^T\g(t)dB_H(t)=\int_{s}^T\g(t)dW_H(t)+\int_{s}^T\g(t)\DD R_H(t)dt\\
=\int_{s}^TG_H(\tau,s,T,\g)d B(\tau)+\int_{s}^T\g(t)\DD R_H(t)dt,
\eaaa
and all integrals here converge  converges  in $L_1(\Omega,\GG_{T},\P)$.
\item[(ii)] $\E\int_{s}^T\g(t)dW_H(t)=0$ but
it may happen that $\E\int_{s}^T\g(t)\DD R_H(t)dt\neq 0$.
\end{enumerate}
\end{corollary}
\begin{remark}
For a small $\Delta s>0$  and $T=s+\Delta s$, the term $W_H(t)_{[s,s+\Delta s]}$ represents a diffusion "noise"
component that is independent on the past, and the term $R_H(t)_{[s,s+\Delta s]}$ represents a
a smooth and predictable component that is completely defined by the past. This and Corollary \ref{corr}9ii) show that
the smooth process $R_H(t)$ has some similarity with
the drift term of a semimartingale  diffusion process.
\end{remark}
In \cite{CNS},\cite{Mu},\cite{Es}, the term "drift" was applied for the parameter $\mu$ in the presentation  $\mu t+B_H(t)$.
\par
In  \cite{Mu}, Theorem 1, a representation of $B_H$ after linear integral transformation and a time change was obtained
via a standard Brownian motion process plus a linear  in time process. Our representation is in a different setting, for the process $B_H$ itself, i.e, without an integral transformation or time change.
\vspace{2mm}
\par
{\em Proof of Theorem \ref{ThM}(i).}
Clearly, the processes $W_H(t)$ and $R_H(t)$  are Gaussian with zero mean and satisfy the required measurability properties.

 By the property of the Riemann--Liouville integral, there exists $c>0$ such that
$$\|G_H(\cdot,s,T,\g)\|_{L_2(s,T)}\le c\|\g\|_{L_2(s,T)}$$ a.s.. In addition, this value
is $\GG_s$-measurable for any $\tau$.
We have  that the integral $\int_{s}^T \g(t)dW_H(\tau)$ is defined as the It\^o's integral
\baa
&&c_H(H-1/2)\int_{s}^T \g(t)dt\int_s^t (t-\tau)^{H-3/2}dB(\tau)\nonumber\\&&=c_H(H-1/2)\int_{s}^T dB(\tau)\int_\tau^T(t-\tau)^{H-3/2}\g(t)dt
=\int_{s}^T dB(\tau)G_H(\tau,s,T,\g)
\label{IW}\eaa
 \par

\index{In particulr, $W_H\to B$ as $H\to 1/2$}

To prove Theorem \ref{ThM}(ii), we need to verify the properties related to the differentiability of
 $R_H(t)$. It suffices to prove the following lemma.
  
\begin{lemma}\label{lemma} The mean square derivative described in Theorem \ref{ThM} exists and can be represented as
 \baaa
 \DD R_H(t)=c_H  \int_{-\infty}^{\Ts}f'_t(t,q)dB(q).
 \eaaa
 For this process,
\baaa
\E\int_{s}^T \DD R_H(t)^2dt= \frac{c_H^2}{4}(T-s)^{2H-1}.\label{EDR}
\eaaa
\end{lemma}

{\em Proof of Lemma \ref{lemma}.}
Let $t>\Ts$ and $q<s$. Consider the derivative
\baaa
f_t'(t,q)=(H-1/2)(t-q)^{H-3/2}.
\eaaa
Since $H-3/2\in (-1,-1/2)$, it follows that  $2(H-3/2)\in (-2,-1)$ and  $\|f'_{t}(t,\cdot)\|_{L_2(-\infty,\Ts)}<+\infty$.

Let $f^{(1)}(t,q,\d)=(f(t+\d,q)-f(t,q))/\d$, where $\d\in (-(t-\Ts)/2,(t-\Ts)/2)$.

Clearly, $f'(t,q)-f^{(1)}(t,q,\d)\to 0$ as $\d\to 0$ for all $t>T_1$ and all $q$.
Let us show that $\|f'(t,\cdot)-f^{(1)}(t,\cdot,\d)\|_{L_2(-\infty,\Ts)}\to 0$ as $\d\to 0$.
We have that
\baaa
f^{(1)}(t,q,\d)=\d^{-1}\int_t^{t+\d}f_t'(s,q)ds=f'(\t(q,\d),q)
\eaaa
for some $\t(q,\d)\in (t,t+\d)$.  Hence
\baa
|f'_t(t,q)-f^{(1)}(t,q,\d)|\le \sup_{h\in (t,t+\d)}|f'_t(t,q)-f'_t(h,q)|\le h\sup_{h\in (t,t+\d)}|f''_{tt}(h,q)|,
\label{est}\eaa
where
\baaa
f''_{tt}(h,q)=(H-1/2)(H-3/2)(h-q)^{H-5/2}.
\eaaa
For $\d>0$, we have that
\baaa
\sup_{h\in (t,t+\d)}|f''_{tt}(h,q)|\le |(H-1/2)(H-3/2)|(t-q)^{H-5/2}.
\eaaa
For $\d\in (-(t-\Ts)/2,0]$, we have that $t+\d-q>t+\d-\Ts>(t-\Ts)/2$, and
\baaa
\sup_{h\in (t,t+\d)}|f''_{tt}(h,q)|\le |(H-1/2)(H-3/2)|(t+\d-q)^{H-5/2}.
\eaaa
Since $H-5/2\in (-2,-3/2)$, it follows that $\|f''_{tt}(t,\cdot)\|_{L_2(-\infty,\Ts)}<+\infty$.

By (\ref{est}), it follows for all $t>s$ that $c_H  \int_{-\infty}^{\Ts}f'_t(t,q)dB(q)$ is mean square limit
\baaa
\lim_{\d\to 0+} \frac{R_H(t+\d)-R_H(t)}{\d}=c_H  \int_{-\infty}^{\Ts}f'_t(t,q)dB(q).
\eaaa
\def\DD{{\cal D}}
We denote this limit as $\DD R_H(t)$, i.e., $\DD R_H(t)=c_H  \int_{-\infty}^{\Ts}f'_t(t,q)dB(q)$.

Further, we have that
\baa
\E \DD R_H(t)^2&=&c_H^2\int_{-\infty}^{\Ts}|f_t'(t,q)|^2dq
=c_H^2(H-1/2)^2\int_{-\infty}^{\Ts}(t-q)^{2H-3}dq\nonumber\\&=&c_H^2\frac{(H-1/2)^2}{2H-2}(t-q)^{2H-2}\Bigl|^{\Ts}_{-\infty}
=c_H^2\frac{H-1/2}{2}(t-\Ts)^{2H-2}.
\label{DR}
\eaa
Hence, for $T>\Ts$,
\baaa
\E\int_{\Ts}^T \DD R_H(t)^2dt =c_H^2\frac{H-1/2}{2}\int_{\Ts}^T (t-\Ts)^{2H-2} dt&=&
c_H^2\frac{H-1/2}{2(2H-1)}(T-\Ts)^{2H-1}\\&=& \frac{c_H^2}{4}(T-\Ts)^{2H-1}.
\eaaa
It follows that (\ref{EDR}) holds.
 This completes the proof of Lemma \ref{lemma} and Theorem \ref{ThM}.  $\Box$
 \par
 {\em Proof of Corollary \ref{corr}}.   The integral $\int_{s}^T \g(t)dW_H(\tau)$ is defined as the It\^o's integral (\ref{IW})
that converges  in $L_2(\Omega,\GG_{T},\P)$. Further,
\baaa
\E\left|\int_{s}^T \g(t)\DD R_H(t)dt\right|\le \left(\E\int_{s}^{T}\g(t)^2dt\right)^{1/2} \left(\E\int_{s}^{T}\DD R_H(t)^2dt\right)^{1/2},
\eaaa and the integral $\int_{s}^T \g(t)\DD R_H(t)dt$ converges  in $L_1(\Omega,\GG_{s},\P)$. Then statement (i) follows. Statement (ii)
follows from a straightforward example where $\g(t)=\DD R_H(t)$. $\Box$
 \section{Applications to financial models}
Consider the market model
consisting of a risk free bond or bank account with the price $b(t)$, ${t\ge 0}$, and
 a risky stock with the price $S(t)$, ${t\ge 0}$. The prices of the stocks evolve
 as \be \label{S} S(t)=S(0)+ \mu t+\s B_H(t),\ee where $B_H(t)$ is a fractional Brownian motion such as described above  with the Hurst exponent $H\in(1/2,1)$. The
initial price $S(0)\in \R$ is given; the parameters and $\mu,\s\in\R$, $\s\neq 0$ are also given.
\par
The price
of the bond evolves as \baaa \label{B} db(t)=r b(t)dt, \eaaa where
$B(0)$ is a  given constant, $r\ge 0$ is a short rate. For simplicity, we assume that $r=0$.
\par
We assume that the wealth $X(t)$ at time $t\in[0,T]$ is
\begin{equation}
\label{X} X(t)=\b(t)b(t)+\g(t)S(t).
\end{equation}
Here $\b(t)$ is the quantity of the bond portfolio, $\g(t)$ is the
quantity of the stock  portfolio, $t\ge 0$. The pair $(\b(\cdot),
\g(\cdot))$ describes the state of the bond-stocks securities
portfolio at time $t$. Each of  these pairs is  called a strategy.

\par  We consider pairs $(\b(\cdot),\g(\cdot))$  such that  $\b(t)$ and $\g(t)$ are progressively
measurable with respect to $\{\GG_t\}$.
In addition, we require that  $$ \E\int_0^{T}\left[\b(t)^2+
\g(t)^2\right]dt<+\infty.$$
This restriction bounds the risk to be accepted and  pays the same role as exclusion of doubling strategies; see examples and
discussion  on doubling strategies in \cite{B11}.

Furthermore, we restrict  consideration by  the set $\oo\A$ consisting of $\g$ such that $\g(t)$ is $\GG_0$-measurable for all $t\ge 0$.
 In other words, we consider
preselected strategies, where decisions are planned ahead at the initial time based on the available
historical data.
This is the setting portfolio selected according to a predetermined program
that cannot be to be adjusted according to the flow of
newly available market information; this helps to exclude arbitrage presented for currently adjustable
strategies
for  market with $H\in (1/2,1)$ (see, e.g., \cite{B11}).
In the stochastic control theory, similar strategies  are called programmed controls.

An  admissible
strategy $(\b(\cdot),\g(\cdot))$  is said to be an admissible
self-financing strategy if
 $$
dX(t)=\b(t)db(t)+\g(t)dS(t)=\g(t)dS(t),
$$
meaning that
 $$
X(t)=X(0)+\int_0^t\b(s)db(s)+\int_0^t\g(s)dS(s),
$$
Since we assumed that $r=0$, it gives simply that
 $$
X(t)=X(0)+\int_0^t\g(s)dS(s)=X(0)+\mu\int_0^t\g(s)ds+\s\int_0^t\g(s) dB_H(s).
$$
Under this condition, the process
$\g(t)$ alone defines the strategy.

The stochastic integral mentioned  above is defined as suggested in Corollary \ref{corr}(i).

Let $\lambda >0$ and $k>0$ be given.
\def\EO{{\bf E}_{\scriptscriptstyle 0}}
 \def\VarO{{\rm Var}_{\scriptscriptstyle 0}}
Let $\E_0$ denote the conditional expectation given $\GG_0$, and
let $\VarO$ denote the corresponding variance.

Consider the following portfolio selection problem
 \baa
\hbox{Maximize}\quad \E X(T)- \lambda \E\VarO X(T)-k\E\int_0^T\g(t)^2dt\quad\hbox{over}\quad \g\in\oo\A.\label{U}
\eaa
The selection of
$\lambda$ is defined by the risk preferences; the selection of
$k>0$ represent a penalty for an excessive investment in the risky asset.

 This linear quadratic optimization setting, or so-called mean-variance portfolio, was
   first introduced for single period models by Markowitz (see \cite{Mar}) and extended later on multi-period and continuous time problems
   (see, e.g., \cite {L2} and \cite{D10} for multi-period problems and \cite{L4} for continuous time).

By Theorem \ref{ThM}, we have that
\baaa
&&\EO X(T)=\int_0^T\g(t)(\mu+\s\DD R_H(t))dt,\qquad \VarO X(T)=\EO\left(\s\int_0^T\g(t) dW_H(t)\right)^2.
\eaaa
Consider an operator $\G:L_2(0,T)\to L_2(0,T)$ such that
\baaa
(\g,\G\g)_{L_2(0,T)}=\s^2(G_H(\cdot,0,T,\g),G_H(\cdot,0,T,\g))_{L_2(0,T)}.
\eaaa
Clearly,
the operator $\G$ is bounded, self-adjoint, and non-negatively defined. By Corollary \ref{corr}, it follows that
\baaa
\E\VarO X(T)=\E(\g,\G\g)_{L_2(0,T)}.
\eaaa
 \par
Consider a Hilbert space $\H=L_2(\O,\GG_0,\P;L_2(0,T))$ with the standard $L_2$-type scalar product $(\cdot,\cdot)_\H$
and norm $\|\cdot\|_\H$.

 Problem (\ref{U}) can be rewritten as
  \baa
\hbox{Maximize}\quad (\g,\mu+\DD R_H)_{\H} - \lambda (\g,\G\g)_{\H}-k \|\g\|^2_{\H}
\quad\hbox{over}\quad \g\in \H.
\label{U2}
\eaa
The  theory of linear-quadratic gives that the optimal solution in $\H$
is
\baa
\w\g= \frac{1}{2\lambda}(\G +k I)^{-1}(\mu + \s\DD R_H).
\label{wg}\eaa
Here $I:L_2(0,T)\to L_2(0,T)$ is the identity operator. The operator  $(\G +k I)^{-1}:L_2(0,T)\to L_2(0,T)$ is bounded
(see, e.g., Lemma 2.7 p.74 in \cite{GT}).
A version of $\w\g$ can be selected such that
 $\g(t)$ is $\GG_0$-measurable for a.e. $t$ and $\g\in L_2(0,T)$ a.s..

\begin{remark}
The operator $(\G +k I)^{-1}$ can be applied to $\mu + \s\DD R_H$ pathwise in $L_2(0,T)$; this does not require
to apply any averaging over a  probability measure.
\end{remark}
\section{Discussion and future developments}
\begin{enumerate}
\item  It can be noted that estimation of $\mu$ and $\DD R_H$ from observable data required by (\ref{wg})   is a non-trivial problem;
see discussion in \cite{CNS}, \cite{Es}, and  \cite{Mu}.
\item It could be interesting to solve a similar portfolio selection problem
for a more mainstream model with $S(t)=\exp(\mu t+\s B_H(t))$. It is unclear yet how to do this for a setting
with $\GG_{0}$-measurable
quantity of shares $\g(t)$ for admissible strategies. However, it is straightforward to consider a model with this prices in a setting with
$\g(t)=\pi(t)/S(t)$, where $\pi(t)$ is $\F_{0}$-adapted. We leave for future research.
\end{enumerate}

\subsection*{Acknowledgment} This work  was
supported by ARC grant of Australia DP120100928 to the author.

\end{document}